\documentclass[12pt]{amsart}

\usepackage{amssymb}
\usepackage{graphicx}
\usepackage{ifthen}
\usepackage{mathrsfs}
\usepackage{pict2e}
\usepackage{xargs}
\usepackage{xspace}

\newcommand{\definedterm}[1]{\emph{#1}}

\newcommand{\Bairespace}[1][]{
  \ifthenelse{\equal{#1}{}}{\functions{\N}{\N}}{\functions{#1}{\N}}
}
\newcommand{\Bairetree}[1][]{
  \ifthenelse{\equal{#1}{}}{\functions{<\N}{\N}}{\functions{#1}{\N}}
}

\newcommand{\boundary}[2]{\partial_{#1}(#2)}

\newcommand{\calB}{\mathcal{B}}
\newcommand{\calG}{\mathcal{G}}

\newcommand{\calX}{\mathcal{X}}
\newcommand{\Cantorspace}[1][]{
  \ifthenelse{\equal{#1}{}}{\functions{\N}{2}}{\functions{#1}{2}}
}
\newcommand{\Cantortree}{\functions{<\N}{2}}
\newcommand{\cardinality}[1]{|#1|}
\newcommandx{\chromaticnumber}[2][1 =]{\chi_{#1}(#2)}
\newcommand{\closedinterval}[2]{[#1, #2]}
\newcommandx{\composition}[2][1 =, 2 =]{
  \ifthenelse{\equal{#1}{}}{\circ}{
    \ifthenelse{\equal{#2}{}}{\circ #1}{{\circ_{#1} #2}}
  }
}
\newcommandx{\concatenation}[2][1 = undefined, 2 = undefined]{
  \ifthenelse{\equal{#1}{undefined}}{{}\smallfrown}{
    \ifthenelse{\equal{#2}{undefined}}{\bigoplus #1}{\bigoplus_{#1} #2}
  }
}
\newcommand{\connectedcomponent}[2]{[#1]_{#2}}
\newcommand{\degree}[2]{\text{deg}_{#1}(#2)}

\newcommand{\emptyfunction}{\emptyset}
\newcommand{\emptysequence}{\emptyset}
\newcommand{\emptystring}{\emptyset}
\newcommand{\equivalenceclass}[2]{[#1]_{#2}}

\newcommand{\extendedby}{\sqsubseteq}

\newcommand{\Ezero}{\mathbb{E}_0}
\newcommand{\forcomeagerlymany}{\forall^*}
\newcommand{\from}{\colon}

\newcommandx{\functions}[3][3 =]{
  \ifthenelse{\equal{#3}{}}{#2^{#1}}{#2^{#1}_{#3}}
}

\newcommand{\generatedgraph}[1]{G_{#1}}

\newcommand{\graph}[1]{\mathrm{graph}(#1)}

\newcommand{\horizontalsection}[2]{#1^{#2}}

\newcommand{\image}[2]{#1(#2)}
\newcommandx{\intersection}[2][1 =, 2 =]{
  \ifthenelse{\equal{#1}{}}{\cap}{
    \ifthenelse{\equal{#2}{}}{\bigcap #1}{{\bigcap_{#1} #2}}
  }
}
\newcommand{\inverse}[1]{#1^{-1}}

\newcommand{\Lzero}{\mathscr{L}_0}
\newcommand{\mathand}{\text{ and }}
\newcommand{\mathcomma}{\text{, }}
\newcommand{\mathcommaand}{\text{, and }}

\newcommand{\N}{\mathbb{N}}

\newcommand{\pair}[2]{(#1, #2)}
\newcommandx{\Piclass}[2][1=,2=]{
  \ifthenelse{\equal{#2}{}}{\mathbf{\Pi}_{#1}}{\mathbf{\Pi}^{#1}_{#2}}
}

\newcommand{\preimage}[2]{#1^{-1}(#2)}
\newcommandx{\projection}[2][1 =, 2 =]{
  \ifthenelse{\equal{#1}{}}{\mathrm{proj}}{
    \ifthenelse{\equal{#2}{}}{\projection_{#1}}{
      \image{\projection[#1]}{#2}
    }
  }
}
\renewcommand{\restriction}[2]{#1 \upharpoonright #2}
\newcommand{\saturation}[2]{[#1]_{#2}}
\newcommandx{\sequence}[2][2 = undefined]{
  \ifthenelse{\equal{#2}{undefined}}{(#1)}{
    (#1)_{#2}
  }
}

\newcommandx{\set}[2][2 = undefined]{
  \ifthenelse{\equal{#2}{undefined}}{\{ #1 \}}{
    \{ #1 \suchthat #2 \}
  }
}

\newcommandx{\sets}[3][3 =]{
  \ifthenelse{\equal{#3}{}}{[#2]^{#1}}{[#2]^{#1}_{#3}}
}
\newcommand{\sextuple}[6]{(#1, #2, #3, #4, #5, #6)}

\newcommandx{\Sigmaclass}[2][1=,2=]{
  \ifthenelse{\equal{#2}{}}{\mathbf{\Sigma}_{#1}}{\mathbf{\Sigma}^{#1}_{#2}}
}

\newcommand{\suchthat}{\mid}

\newcommand{\tailequivalence}[1]{E_t(#1)}

\newcommand{\transposition}[2]{(#1 \ #2)}

\newcommand{\Tzero}{\mathscr{T}_0}
\newcommandx{\union}[2][1 =, 2 =]{
  \ifthenelse{\equal{#1}{}}{\cup}{
    \ifthenelse{\equal{#2}{}}{\bigcup #1}{{\bigcup_{#1} #2}}
  }
}
\newcommand{\verticalsection}[2]{#1_{#2}}

\newcommand{\Adams}{Ad\-ams\xspace}
\newcommand{\Baire}{Baire\xspace}
\newcommand{\Bernoulli}{Ber\-noul\-li\xspace}
\newcommand{\Borel}{Bor\-el\xspace}
\newcommand{\Cayley}{Cay\-ley\xspace}
\newcommand{\Dougherty}{Dough\-er\-ty\xspace}
\newcommand{\Harrington}{Har\-ring\-ton\xspace}
\newcommand{\Hjorth}{Hjorth\xspace}
\newcommand{\Jackson}{Jack\-son\xspace}
\newcommand{\Kechris}{Kech\-ris\xspace}

\newcommand{\Kuratowski}{Kur\-at\-ow\-ski\xspace}
\newcommand{\Laczkovich}{Lacz\-ko\-vich\xspace}
\newcommand{\Lebesgue}{Leb\-esgue\xspace}
\newcommand{\Louveau}{Lou\-veau\xspace}
\newcommand{\Lusin}{Lu\-sin\xspace}
\newcommand{\Lyons}{\-Ly\-ons\xspace}
\newcommand{\Marks}{Marks\xspace}

\newcommand{\Nazarov}{Naz\-ar\-ov\xspace}
\newcommand{\Novikov}{No\-vik\-ov\xspace}
\newcommand{\Polish}{Po\-lish\xspace}
\newcommand{\Slaman}{Sla\-man\xspace}
\newcommand{\Steel}{Steel\xspace}
\newcommand{\Ulam}{U\-lam\xspace}
\newcommand{\Weiss}{Weiss\xspace}

\newenvironment{lemmaproof}{
  
  \begin{proof}
}{\end{proof}}

\newenvironment{propositionproof}{
  
  \begin{proof}
}{\end{proof}}

\newenvironment{theoremproof}{
  
  \begin{proof}
}{\end{proof}}

\newtheorem{corollary}{Corollary}[section]
\newtheorem{lemma}[corollary]{Lemma}
\newtheorem{proposition}[corollary]{Proposition}
\newtheorem{theorem}[corollary]{Theorem}

\newtheorem{introconjecture}{Conjecture}
\newtheorem{introtheorem}[introconjecture]{Theorem}

\theoremstyle{definition}

\newtheorem{remark}[corollary]{Remark}

\begin{document}


\begin{abstract}
  We characterize the structural impediments to the existence of \Borel perfect
  matchings for acyclic locally countable \Borel graphs admitting a \Borel
  selection of finitely many ends from their connected components. In particular,
  this yields the existence of \Borel matchings for such graphs of degree at least
  three. As a corollary, it follows that acyclic locally countable \Borel graphs of
  degree at least three generating $\mu$-hyperfinite equivalence relations admit
  $\mu$-measurable matchings. We establish the analogous result for \Baire
  measurable matchings in the locally finite case, and provide a counterexample
  in the locally countable case.
\end{abstract}

\author[C.T. Conley]{Clinton T. Conley}

\address{
  Clinton T. Conley \\
  Department of Mathematical Sciences \\
  Carnegie Mellon University \\
  Pittsburgh, PA 15213 \\
  USA
}

\email{clintonc@andrew.cmu.edu}

\urladdr{
  http://www.math.cmu.edu/math/faculty/Conley
}

\author[B.D. Miller]{Benjamin D. Miller}

\address{
  Benjamin D. Miller \\
  Kurt G\"{o}del Research Center for Mathematical Logic \\
  Universit\"{a}t Wien \\
  W\"{a}hringer Stra{\ss}e 25 \\
  1090 Wien \\
  Austria
 }

\email{benjamin.miller@univie.ac.at}

\urladdr{
  http://www.logic.univie.ac.at/benjamin.miller
}

\thanks{The first author was support in part by NSF Grant DMS-1500906.}

\thanks{The second author was supported in part by DFG SFB Grant 878 and FWF
  Grant P28153.}



\title[Measurable perfect matchings]{Measurable perfect matchings for acyclic
  locally countable Borel graphs}

\maketitle

\section*{Introduction} \label{introduction}

A \definedterm{graph} on a set $X$ is an irreflexive symmetric subset $G$ of
$X \times X$. An \definedterm{involution} is a permutation which is its own
inverse, and a \definedterm{matching} of $G$ is an involution of a subset of
$X$ whose graph is contained in $G$. Such a matching is \definedterm{perfect}
if its domain is $X$ itself.

A \definedterm{$G$-path} is a sequence $\sequence{x_i}[i \le n]$ such that $x_i
\mathrel{G} x_{i+1}$, for all $i < n$. We say that $G$ is \definedterm{connected}
if there is a $G$-path between any two points of $X$. More generally, the equivalence relation
\definedterm{generated} by a graph $G$ on $X$ is the smallest equivalence relation
on $X$ containing $G$, and the \definedterm{connected components} of $G$ are the
equivalence classes $\connectedcomponent{x}{G}$ of this relation. A graph $G$ is 
\definedterm{acyclic} if there is at most one injective $G$-path between any two points.
When $G$ is acyclic, the \definedterm{$G$-distance} between two points of the same
connected component of $G$ is one less than the number of points along the unique 
injective $G$-path between them. A \definedterm{tree} is an acyclic connected graph.

A straightforward recursive analysis yields a characterization of the existence of perfect
matchings for acyclic graphs. Here we consider the substantially more subtle question
of the existence of measurable perfect matchings for acyclic definable graphs.

A \definedterm{\Polish space} is a separable topological space admitting a
compatible complete metric.  A subset of such a space is \definedterm{\Borel} if
it is in the $\sigma$-algebra generated by the underlying topology. A \definedterm
{standard \Borel space} is a set $X$ equipped with the family of \Borel sets
associated with a \Polish topology on $X$. Every subset of a standard \Borel
space inherits the $\sigma$-algebra consisting of its intersection with each
\Borel subset of the original space; this restriction is again standard \Borel exactly
when the subset in question is \Borel (see, for example, \cite[Corollary 13.4 and
Theorem 15.1]{Kechris}). A function between standard \Borel spaces is
\definedterm{\Borel} if pre-images of \Borel sets are \Borel. We will take being 
\Borel as our notion of definability.

The \definedterm{$G$-degree} of a point $y$ is given by $\degree{G}{y} =
\cardinality{\set{x \in X}[x \mathrel{G} y]}$. A graph is \definedterm{locally
countable} if every point has countable $G$-degree, and \definedterm{locally
finite} if every point has finite $G$-degree. A graph is \definedterm{$n$-regular} if
every point has $G$-degree $n$. We say that a graph has \definedterm{degree at
least $n$} if every point has $G$-degree at least $n$. The existence of perfect
matchings can be reduced to the case of graphs of degree at least two (modulo a
minor caveat in the \Borel setting).

A \definedterm{$G$-ray} is a sequence $\sequence{x_n}[n \in \N]$ with the property
that $x_n \mathrel{G} x_{n+1}$, for all $n \in \N$. We say that a sequence $\sequence
{x_n}[n \in \N]$ has \definedterm{$G$-degree two on even indices} if $\degree{G}
{x_{2n}} = 2$, for all $n \in \N$. Note that if $G$ has degree at least three, then
there are no such sequences.

When $G$ is acyclic, we say that injective $G$-rays $\sequence{x_n}[n
\in \N]$ and $\sequence{y_n}[n \in \N]$ are \definedterm{end equivalent} if
there exist $i, j \in \N$ with $x_{i + n} = y_{j + n}$, for all $n \in \N$. We say that
a set $\calX \subseteq \functions{\N}{X}$ \definedterm{selects} a finite non-empty
set of ends from every connected component of $G$ if $\calX \intersection \functions
{\N}{\connectedcomponent{x}{G}}$ is a finite non-empty union of end-equivalence
classes, for all $x \in X$.

\begin{introtheorem} \label{introduction:Borel}
  Suppose that $X$ is a \Polish space, $G$ is an acyclic locally countable \Borel
  graph on $X$ of degree at least two, and there is a \Borel set selecting a finite
  non-empty set of ends from every connected component of $G$. Then there is
  a \Borel set $B \subseteq X$ such that:
  \begin{enumerate}
    \item The restriction of $G$ to $B$ is two-regular.
    \item Every connected component of $G$ contains at most one connected
      component of the restriction of $G$ to $B$.
    \item No two points of $B$ of $G$-degree at least three have odd $G$-distance
      from one another.
    \item There is a \Borel perfect matching of $G$ off of $B$.
  \end{enumerate}
  In particular, if there are no injective $G$-rays of $G$-degree two on even indices,
  then the set $B$ is empty, thus $G$ has a \Borel perfect matching.
\end{introtheorem}

There are well-known examples of acyclic two-regular \Borel graphs which do not have
\Borel perfect matchings, and a result of \Marks yields acyclic $n$-regular \Borel
graphs which do not have \Borel perfect matchings, for all natural numbers $n \ge
3$ (see \cite[Theorem 1.5]{Marks}).

A \definedterm{\Borel probability measure} on a \Polish space $X$ is a
function $\mu$, assigning to each \Borel set $B \subseteq X$ an element of
$\closedinterval{0}{1}$, with the property that $\mu(X) = 1$ and $\mu(\union
[n \in \N][B_n]) = \sum_{n \in \N} \mu(B_n)$, for every sequence $\sequence
{B_n}[n \in \N]$ of pairwise disjoint \Borel subsets of $X$. A \Borel set $B
\subseteq X$ is \definedterm{$\mu$-null} if $\mu(B) = 0$, and \definedterm
{$\mu$-conull} if its complement is $\mu$-null.

Following the usual abuse of language, we say that an equivalence relation is
\definedterm{countable} if its classes are countable, and \definedterm{finite} if its
classes are finite. We say that a countable \Borel equivalence relation on a standard
\Borel space is \definedterm{hyperfinite} if it is the union of an increasing sequence
$\sequence{F_n}[n \in \N]$ of finite \Borel subequivalence relations, and
\definedterm{$\mu$-hyperfinite} if there is a $\mu$-conull invariant \Borel set on
which it is hyperfinite.

A well-known result of \Adams and \Jackson-\Kechris-\Louveau (see \cite[Lemma 3.21]
{JacksonKechrisLouveau}) ensures that if $G$ is an acyclic locally countable \Borel
graph on $X$, then the equivalence relation generated by $G$ is $\mu$-hyperfinite if and
only if there is a $\mu$-conull $G$-invariant \Borel set on which there is a \Borel set 
selecting a finite non-empty set of ends from every connected component of $G$ that
has injective $G$-rays. Theorem \ref{introduction:Borel} therefore yields the following
corollary.

\begin{introtheorem} \label{introduction:measure}
  Suppose that $X$ is a \Polish space, $G$ is an acyclic locally countable \Borel
  graph on $X$ of degree at least two and with no injective $G$-rays of $G$-degree
  two on even indices, and $\mu$ is a \Borel probability measure on $X$ for which
  the equivalence relation generated by $G$ is $\mu$-hyperfinite. Then there is
  a $\mu$-conull $G$-invariant \Borel set on which $G$ has a \Borel perfect
  matching.
\end{introtheorem}

By a result of \Lyons-\Nazarov, a wide class of regular bipartite \Borel graphs
(notably including any bipartite \Cayley graphing of the \Bernoulli shift action of a
nonamenable group) admit $\mu$-measurable matchings (see \cite{LyonsNazarov}).
The general case of $\mu$-measurable matchings for (not necessarily bipartite)
\Cayley graphings of \Bernoulli shifts of nonamenable groups is discussed in
\cite{CsokaLippner}.

A subset of a \Polish space is \definedterm{meager} if it is a countable union of
nowhere dense sets, and \definedterm{comeager} if its complement is meager.
In contrast with the measure-theoretic setting, a well-known result of \Hjorth-\Kechris
implies that every countable \Borel equivalence relation is hyperfinite on a comeager
invariant \Borel set (see, for example, \cite[Theorem 12.1]{KechrisMiller}). However,
there are acyclic locally finite \Borel graphs of degree at least two which do not admit
\Borel sets selecting a finite non-empty set of ends on any comeager invariant
\Borel set (see, for example, the graph $\Tzero$ of \cite{HjorthMiller}).
Nevertheless, an entirely different approach yields an analog of Theorem \ref
{introduction:measure} in this context.

\begin{introtheorem} \label{introduction:category}
  Suppose that $X$ is a \Polish space and $G$ is an acyclic locally finite \Borel
  graph on $X$ of degree at least two and with no injective $G$-rays of $G$-degree
  two on even indices. Then there is a comeager $G$-invariant \Borel set on which
  $G$ has a \Borel perfect matching.
\end{introtheorem}

However, we provide an example of an $\aleph_0$-regular \Borel graph 
which does not have a \Borel perfect matching on a comeager invariant \Borel set.
Some rather general sufficient conditions for the existence of \Baire measurable
matchings of graphs are presented in \cite{KechrisMarks} and \cite{MarksUnger}.

The paper is organized as follows. In \S\ref{choice}, we mention a pair of
elementary facts concerning matchings outside of the definable context. In \S
\ref{Borel}, we establish Theorem \ref{introduction:Borel}. And in \S\ref{measurable},
we establish Theorems \ref{introduction:measure} and Theorem \ref{introduction:category},
and describe the example mentioned above.

\section{Matchings using choice} \label{choice}

Here we establish a pair of elementary facts, whose proofs will later prove useful
in the definable setting.

Clearly the existence of a perfect matching for a graph on a non-empty set
necessitates that the graph in question has degree at least one. The following
observation allows one to focus upon graphs of degree at least two.

\begin{proposition} \label{choice:derivative}
  Suppose that $X$ is a set and $G$ is a graph on $X$. Then there is a set $Y \subseteq
  X$, on which $G$ has degree at least two, with the property that if there is a matching
  $\iota$ of $G$ whose domain contains $X \setminus Y$, then $X \setminus Y$ is an
  $\iota$-invariant set on which every other such matching agrees with $\iota$.
\end{proposition}

\begin{propositionproof}
  The \emph{$G$-boundary} of a set $Y \subseteq X$, denoted by $\boundary{G}{Y}$,
  is the set of points in $Y$ which are $G$-related to at least one point outside of $Y$.
  Let $\alpha$ denote the supremum of the ordinals $\beta$ for which there exists $x
  \in X$ such that $\cardinality{\beta} \le \cardinality{\connectedcomponent{x}{G}}$,
  and recursively define a decreasing sequence $\sequence{X^\beta}[\beta \le \alpha]$
  of subsets of $X$ by setting $X^0 = X$, $X^\lambda = \intersection[\beta < \lambda]
  [X^\beta]$, $X^{\lambda + 2n + 1} = \set{x \in X^{\lambda + 2n}}[\degree{\restriction
  {G}{X^{\lambda + 2n}}}{x} \ge 2]$, and $X^{\lambda + 2n + 2} = X^{\lambda + 2n + 1}
  \setminus \boundary{\restriction{G}{X^{\lambda + 2n}}}{X^{\lambda + 2n} \setminus
  X^{\lambda + 2n + 1}}$, for all limit ordinals $\lambda$ and natural numbers $n$ such
  that the corresponding indices are at most $\alpha$.
  
  Set $Y = X^\alpha$. As $X^\alpha = X^{\alpha + 1}$, it follows that $\restriction{G}{Y}$
  has degree at least two. And a straightforward transfinite induction shows that if $\iota$ is a
  matching of $G$ whose domain contains $X \setminus Y$, then $\iota(x)$ is the unique
  $G$-neighbor of $x$ in $X^{\lambda + 2n}$ for all limit ordinals $\lambda$,
  natural numbers $n$, and $x \in X^{\lambda + 2n} \setminus X^{\lambda + 2n + 1}$,
  whereas $\iota(x)$ is the unique $G$-neighbor of $x$ in $X^{\lambda + 2n}
  \setminus X^{\lambda + 2n + 1}$ for all limit ordinals $\lambda$, natural numbers
  $n$, and $x \in X^{\lambda + 2n + 1} \setminus X^{\lambda + 2n + 2}$.
\end{propositionproof}

In the absence of definability requirements, the following completes the analysis of the
existence of perfect matchings for acyclic graphs.

\begin{proposition} \label{choice:matching}
  Suppose that $X$ is a set and $G$ is an acyclic graph on $X$ of degree at least one
  whose connected components have at most one point of $G$-degree exactly one.
  Then $G$ has a perfect matching. 
\end{proposition}

\begin{propositionproof}
  A \definedterm{transversal} of an equivalence relation is a set intersecting every
  equivalence class in exactly one point. We will recursively define a sequence
  $\sequence{X_n}[n \in \N]$ of pairwise disjoint subsets of $X$, as well as a sequence
  $\phi_n \from X_{2n} \to X_{2n+1}$ of functions whose graphs are contained in $G$,
  with the property that the graph $G_n = \restriction{G}{(X \setminus \union[m < 2n]
  [X_m])}$ has degree at least one and $X_{2n}$ is a transversal of the
  equivalence relation generated by $G_n$, containing every point of $G_n$-degree
  exactly one.

  We begin by fixing a transversal $X_0 \subseteq X$ of the equivalence relation
  generated by $G$, containing every point of $G$-degree exactly one. Suppose
  now that $n \in \N$ and we have already found $\sequence{X_m}[m \le 2n]$. Fix
  a function $\phi_n \from X_{2n} \to X \setminus \union[m \le 2n][X_m]$
  whose graph is contained in $G$, and set $X_{2n+1} = \image{\phi_n}{X_{2n}}$
  and $X_{2n+2} = \boundary{G}{X \setminus \union[m \le 2n + 1][X_m]}$.
  
  Set $\phi = \union[n \in \N][\phi_n]$, and observe that the involution $\iota = \phi
  \union \inverse{\phi}$ is a perfect matching of $G$.
\end{propositionproof}

\section{Borel matchings} \label{Borel}

We will frequently employ the following well-known fact.

\begin{theorem}[\Lusin-\Novikov] \label{Borel:uniformization}
  Suppose that $X$ and $Y$ are \Polish spaces and $R \subseteq X \times Y$
  is a \Borel set whose vertical sections are all countable. Then $\image{\projection
  [X]}{R}$ is \Borel and there are \Borel functions $\phi_n \from \image{\projection[X]}
  {R} \to Y$ such that $R = \union[n \in \N][\graph{\phi_n}]$.
\end{theorem}

\begin{theoremproof}
  See, for example, \cite[Theorem 18.10]{Kechris}.
\end{theoremproof}

There is a natural analog of Proposition \ref{choice:derivative} in the \Borel setting.

\begin{proposition} \label{Borel:derivative}
  Suppose that $X$ is a \Polish space and $G$ is a locally finite \Borel graph on $X$. Then
  there is a \Borel set $B \subseteq X$, on which $G$ has degree at least two, with
  the property that if there is a matching $\iota$ of $G$ whose domain contains
  $X \setminus B$, then the latter is an $\iota$-invariant set on which every other
  such matching agrees with $\iota$. In particular, it follows that the restriction of
  every such matching to $X \setminus B$ is \Borel.
\end{proposition}

\begin{propositionproof}
  Following the proof of Proposition \ref{choice:derivative}, our assumption that $G$ is
  locally finite ensures that $X_\alpha = X_\omega$, and Theorem \ref{Borel:uniformization}
  implies that the set $B = X_\omega$ is \Borel.
  
  Again by Theorem \ref{Borel:uniformization}, there are \Borel functions $\phi_n \from X
  \to X$ such that the equivalence relation generated by $G$ is $\union[n \in \N][\graph
  {\phi_n}]$. We say that a function $i \from \N \to \N$ \definedterm{codes} a function on
  the connected component of $x$ off of $B$ if $\phi_{i(m)}(x) = \phi_{i(n)}(x)$ and
  neither is in $B$, for all $m, n \in \N$ such that $\phi_m(x) = \phi_n(x)$ and neither is
  in $B$. The corresponding function $\iota \from \connectedcomponent{x}{G} \setminus
  B \to \connectedcomponent{x}{G} \setminus B$ is then given by
  \begin{equation*}
    \iota(y) = z \iff \exists m, n \in \N \ (y = \phi_m(x) \mathcomma z = \phi_n(x)
      \mathcommaand i(m) = n).
  \end{equation*}
  
  Note that if $\iota$ is a matching of $G$ on $X \setminus B$, then $\iota(y) = z$ if
  and only if there is a function $i \from \N \to \N$, coding a matching of $G$ on the
  connected component of $x$ off of $B$, which sends $y$ to $z$. As the set of
  $\pair{x}{i} \in X \times I$ for which $i$ codes a matching of $G$ on the
  connected component of $x$ off of $B$ is \Borel, as is the set of $\sextuple
  {i}{m}{n}{x}{y}{z} \in \Bairespace \times \N \times \N \times X \times X \times X$ for
  which $y = \phi_m(x)$, $z = \phi_n(x)$, and $i(m) = n$, it follows that the graph of
  $\iota$ is \definedterm{analytic}, in the sense that it is the image of a \Borel subset
  of a standard \Borel space under a \Borel function. As functions between \Polish
  spaces with analytic graphs are \Borel (see, for example, \cite[Theorem 14.12]
  {Kechris}), it follows that $\iota$ is \Borel.
\end{propositionproof}

\begin{remark}
  Proposition \ref{choice:derivative} also has an analog in the more general setting
  of locally countable \Borel graphs, granting that we slightly relax the requirement that
  the matching is \Borel. To be specific, in this case one can check that the set $A =
  X_{\alpha}$ is analytic, and that the graph of the unique matching of $\restriction
  {G}{(X \setminus A)}$ is both relatively analytic and co-analytic. Here it is worth
  noting that our other results generalize to \Borel graphs on analytic sets. It is
  also worth noting that, off of a meager or $\mu$-null $G$-invariant \Borel set,
  this yields the full conclusion of Proposition \ref{Borel:derivative}.
\end{remark}

Proposition \ref{choice:matching} also has a natural analog in the \Borel setting,
albeit only when the equivalence relation generated by the graph in question is
particularly simple. A \definedterm{reduction} of an equivalence relation $E$ on $X$ to an
equivalence relation $F$ on $Y$ is a function $\pi \from X \to Y$ such that $x_1 \mathrel{E}
x_2 \iff \pi(x_1) \mathrel{F} \pi(x_2)$, for all $x_1, x_2 \in X$. A \Borel equivalence relation
on a \Polish space is \definedterm{smooth} if it is \Borel reducible to equality on a \Polish
space.

\begin{proposition} \label{Borel:smoothmatching}
  Suppose that $X$ is a \Polish space and $G$ is an acyclic locally countable \Borel graph
  on $X$ of degree at least one, each of whose connected components have at most one
  point of $G$-degree exactly one, whose induced equivalence relation is smooth.
  Then $G$ has a \Borel perfect matching.
\end{proposition}

\begin{propositionproof}
  Following the proof of Proposition \ref{choice:matching}, Theorem \ref{Borel:uniformization}
  ensures that we can choose the functions $\phi_n$ and the sets $X_n$ to be \Borel, in
  which case the corresponding matching is also \Borel.
\end{propositionproof}

Beyond the smooth case, the purely combinatorial and definable settings are
quite different. The graph \definedterm{generated} by a function $f \from X \to X$ is
the graph $\generatedgraph{f}$ on $X$ with respect to which two distinct points
are related if $f$ sends one to the other. A function $f \from X \to X$ is \definedterm
{aperiodic} if $f^n$ is fixed-point free, for all $n > 0$.

\begin{proposition}[\Laczkovich]
  There is a \Polish space $X$ and an aperiodic \Borel automorphism $T \from X \to X$
  with the property that $\generatedgraph{T}$ does not have a \Borel perfect matching.
\end{proposition}

\begin{propositionproof}
  We say that a \Borel probability measure $\mu$ on $X$ is \definedterm{$T$-quasi-invariant}
  if $\mu(B) = 0 \iff \mu(\image{T}{B}) = 0$, for all \Borel sets $B \subseteq X$. And we
  say that a \Borel probability measure $\mu$ on $X$ is \definedterm{$T$-ergodic} if $\mu(B)
  \in \set{0, 1}$, for all $T$-invariant \Borel sets $B \subseteq X$.

  An \definedterm{$I$-coloring} of a graph $G$ on $X$ is a function $c \from X \to I$ with
  the property that $\forall \pair{x}{y} \in G \ c(x) \neq c(y)$. It is easy to see that
  $\generatedgraph{T}$ has a \Borel perfect matching if and only if $\generatedgraph{T}$
  has a \Borel two-coloring, or equivalently, if there is a \Borel set $B \subseteq X$ such that
  $B$ and $\image{T}{B}$ partition $X$. And the latter is ruled out by the existence of a
  \Borel probability measure on $X$ which is $T$-quasi-invariant and $T^2$-ergodic.
    
  As \Lebesgue measure is well known to be ergodic and quasi-invariant with respect to
  irrational rotations of the circle, it follows that the latter do not have \Borel (or even
  \Lebesgue measurable) matchings.
\end{propositionproof}

\begin{remark}
  The above argument goes through just as well using \Baire category in lieu of \Lebesgue
  measure.
\end{remark}

\begin{remark}
  The existence of such measures (or topologies with corresponding notions of \Baire
  category) for graphs generated by aperiodic \Borel automorphisms is, in fact, equivalent
  to the inexistence of \Borel perfect matchings. This follows from a dichotomy theorem
  of \Louveau's (see, for example, \cite[Theorem 15]{Miller:Survey}).
\end{remark}

Let $\Ezero$ denote the equivalence relation on $\Cantorspace$ given by
\begin{equation*}
  x \mathrel{\Ezero} y \iff \exists n \in \N \forall m \ge n \ x(m) = y(m).
\end{equation*}
The \Harrington-\Kechris-\Louveau $\Ezero$ dichotomy ensures that, under \Borel
reducibility, this is the minimal non-smooth \Borel equivalence relation (see \cite
[Theorem 1.1]{HarringtonKechrisLouveau}). Arguments of \Dougherty-\Jackson-\Kechris
can be used to show that a countable \Borel equivalence relation on a \Polish space is
hyperfinite if and only if it is \Borel reducible to $\Ezero$ (see, for example, \cite[Theorem
1]{DoughertyJacksonKechris}), and \Slaman-\Steel and \Weiss have noted that every \Borel
automorphism of a \Polish space generates a hyperfinite \Borel equivalence relation (see,
for example, \cite[Theorem 5.1]{DoughertyJacksonKechris}). In particular, it follows that
the smoothness of the equivalence relation in Proposition \ref{Borel:smoothmatching}
cannot be weakened.

Among graphs generated by aperiodic \Borel functions, there are essentially no further
examples without \Borel matchings. The \definedterm{tail equivalence relation} on $X$
induced by a function $f \from X \to X$ is given by
\begin{equation*}
  x \mathrel{\tailequivalence{f}} y \iff \exists m, n \in \N \ f^m(x) = f^n(y).
\end{equation*}
Note that if $G$ is the graph generated by $f$, then $\tailequivalence{f}$ is the
equivalence relation generated by $G$.

The \definedterm{injective part} of $f$ is the set $\set{x \in X}[\restriction{f}
{\equivalenceclass{x}{\tailequivalence{f}}} \text{ is injective}]$. When $f$ is
countable-to-one, Theorem \ref{Borel:uniformization} ensures that the injective
part of $f$ is a \Borel set on which $f$ is a \Borel automorphism.

\begin{proposition} \label{Borel:function}
  Suppose that $X$ is a \Polish space and $f \from X \to X$ is a countable-to-one \Borel
  surjection. Then $\generatedgraph{f}$ has a \Borel perfect matching, off of the
  injective part of $f$.
\end{proposition}

\begin{propositionproof}
  By Theorem \ref{Borel:uniformization}, there is a \Borel function $g \from X \to X$
  such that $(f \composition g)(x) = x$, for all $x \in X$. Off of the $\tailequivalence
  {f}$-saturation of $\intersection[n \in \N][\image{g^n}{X}]$, the involution agreeing
  with $g$ on $\union[n \in \N][\image{g^{2n}}{X \setminus \image{g}{X}}]$ and with
  $f$ on $\union[n \in \N][\image{g^{2n+1}}{X \setminus \image{g}{X}}]$ is a \Borel
  perfect matching of $\generatedgraph{f}$. So it only remains to produce a \Borel
  perfect matching of $\generatedgraph{f}$ on the $\tailequivalence{f}$-saturation
  of the set
  \begin{equation*}
    \textstyle
    B = \set{x \in \intersection[n \in \N][\image{g^n}{X}]}[x \text{ is not in the injective
      part of } f].
  \end{equation*}
  
  Towards this end, define $A = \set{x \in B}[\cardinality{\preimage{f}{x}} \ge 2]$. As Proposition 
  \ref{Borel:smoothmatching} allows us to throw out an $\tailequivalence{f}$-invariant \Borel set
  on which $\tailequivalence{f}$ is smooth, we can assume that for all $x \in B$, there exist $m,
  n \in \N$ such that $(\restriction{f}{B})^{-m}(x), f^n(x) \in A$. For each $x \in A$, let $n(x)$ denote
  the least positive natural number $n$ for which $f^n(x) \in A$, and define $A' = \set{x \in A}[n(x) 
  \text{ is odd}]$. We then obtain a \Borel perfect matching of $\generatedgraph{f}$ on $B
  \setminus A'$ by associating $f^{2i+1}(x)$ with $f^{2i+2}(x)$ for all $i \in \N$ and $x \in A'$ with
  $2i+2 < n(x)$, as well as $f^{2i}(x)$ with $f^{2i+1}(x)$ for all $i \in \N$ and $x \in A \setminus A'$
  with $2i+1 < n(x)$. As the equivalence relation generated by the restriction of $\generatedgraph
  {f}$ to $A' \union (\saturation{B}{\tailequivalence{f}} \setminus B)$ is smooth, Proposition \ref
  {Borel:smoothmatching} yields an extension to a \Borel perfect matching of $\generatedgraph
  {f}$ on $\saturation{B}{\tailequivalence{f}}$.
\end{propositionproof}

We now turn our attention to another class of \Borel graphs without \Borel perfect matchings.
The \definedterm{line-and-point graph} associated with a graph $G$ on $X$ is the graph on
the disjoint union of $X$ with the set $E = \set{\set{x, y}}[x \mathrel{G} y]$ of \definedterm
{unordered edges of G}, in which two elements of $X \union E$ are related if one of them is
in $X$, one of them is in $E$, and the former is an element of the latter.

Note that if $G$ is a \Borel graph on a standard \Borel space, then the set of unordered edges
of $G$ inherits a standard \Borel structure from $G$, thus the line-and-point graph of $G$ can
also be viewed as a \Borel graph on a standard \Borel space.

Observe also that if a graph has an injective ray, then its line-and-point graph has an injective ray
of degree two on even indices. Together with the following proposition, this is part of the
motivation for focusing on graphs without such rays in our later results.

\begin{proposition} \label{Borel:lineandpoint}
  Suppose that $X$ is a \Polish space and $G$ is an acyclic \Borel graph on $X$. Then $G$
  is generated by an aperiodic \Borel function if and only if its line-and-point graph has a 
  \Borel perfect matching.
\end{proposition}

\begin{propositionproof}
  If $f \from X \to X$ is an aperiodic function generating $G$, then the fact that $f$ is
  fixed-point free ensures that $\set{x, f(x)}$ is an unordered edge of $G$ for all $x \in X$,
  and the fact that $f^2$ is fixed-point free ensures that the involution $\iota$ associating
  $x$ with $\set{x, f(x)}$ is injective. As the fact that $f$ generates $G$ ensures that $\iota$
  is surjective, it is necessarily a perfect matching of the line-and-point graph of $G$.
  
  Conversely, if $\iota$ is a perfect matching of the line-and-point graph of $G$, then the
  function $f$, sending each point $x$ to the unique point $y$ with the property that $\iota
  (x) = \set{x, y}$, generates $G$. The definition of $f$ ensures that both $f$ and $f^2$
  are fixed-point free, and the acyclicity of $G$ ensures that $f^n$ is fixed-point free for
  all $n > 2$, thus $f$ is aperiodic.
\end{propositionproof}

\begin{remark}
  One can drop the acyclicity of $G$ in the statement of Proposition \ref{Borel:lineandpoint}
  by weakening the hypothesis that $G$ is generated by an aperiodic \Borel function to the 
  hypothesis that $G$ is generated by a \Borel function for which both $f$ and $f^2$ are
  fixed-point free.
\end{remark}

\begin{remark}
  When $G$ is an acyclic locally countable \Borel graph of degree at least two, the
  hypothesis that $G$ is generated by an aperiodic \Borel function is equivalent to the
  apparently weaker hypothesis that $G$ is generated by a \Borel function.
\end{remark}

We next consider the combinatorially simplest examples of \Borel graphs which are not induced
by \Borel functions.

\begin{proposition}
  There is a \Polish space $X$ and an acyclic two-regular \Borel graph $G$ on $X$ which is
  not induced by a \Borel function, thus there is such a graph which does not have a \Borel
  perfect matching.
\end{proposition}

\begin{propositionproof}
  The graph $\Lzero$ of \cite{HjorthMiller} yields an example of an acyclic two-regular \Borel
  graph on a \Polish space which is not induced by a \Borel function. Proposition \ref
  {Borel:lineandpoint} then ensures that the corresponding line-and-point graph has no
  \Borel perfect matching (and clearly it is not induced by a \Borel function, since the restriction
  of the square of an aperiodic such function to $2 \times \Cantorspace$ would generate $\Lzero$).
\end{propositionproof}

\begin{remark}
  It is not difficult to verify that the fact we used in the parenthetical remark above is far more
  general. Namely, a \Borel graph on a \Polish space is generated by a \Borel function if
  and only if its line-and-point graph is generated by a \Borel function.
\end{remark}

Among graphs for which such combinatorially simple graphs can be isolated in
a \Borel fashion, there are essentially no further examples without \Borel matchings.

\begin{proposition} \label{Borel:line}
  Suppose that $X$ is a \Polish space, $G$ is an acyclic locally countable \Borel graph on $X$
  of degree at least two, and there is a \Borel set $B \subseteq X$ such that:
  \begin{enumerate}
    \item The restriction of $G$ to $B$ is two-regular.
    \item Every connected component of $G$ contains exactly one connected component of
      the restriction of $G$ to $B$.
    \item Every connected component of $G$ contains two points in $B$, with $G$-degree at
      least three, having odd $G$-distance from one another.
  \end{enumerate}
  Then there is a \Borel perfect matching of $G$.
\end{proposition}

\begin{propositionproof}
  Let $A$ denote the set of points of $B$ of $G$-degree at least three, and let $A'$ denote the
  set of initial points of injective $G$-paths whose initial and terminal points are in $A$, whose other 
  points are not in $A$, and along which there are an even number of points. As Proposition \ref
  {Borel:smoothmatching} allows us to throw out a $G$-invariant \Borel set on which the equivalence 
  relation generated by $G$ is smooth, we can assume that for all $x \in B$, there are points of $A'$ on 
  either \definedterm{side} of $x$, in the sense that for both $G$-neighbors $y$ of $x$ in $B$, there is
  an injective $G$-path of the form $\sequence{x, y, \ldots}$ whose terminal point is in $A'$.
  By Theorem \ref{Borel:uniformization}, there is a \Borel set $A'' \subseteq A'$ consisting of exactly
  one point from every pair of points in $A'$ between which there is an injective $G$-path whose
  other points are not in $A'$ and along which there are an odd number of points. Then there is a
  \Borel perfect matching of the restriction of $G$ to $A'' \union (B \setminus A')$, and Proposition 
  \ref{Borel:smoothmatching} yields an extension of the latter to a \Borel perfect matching of $G$.
\end{propositionproof}

We can now establish the main result of this section.

\begin{theorem} \label{Borel:ends}
  Suppose that $X$ is a \Polish space, $G$ is an acyclic locally countable \Borel
  graph on $X$ of degree at least two, and there is a \Borel set selecting a finite
  non-empty set of ends from every connected component of $G$. Then there is
  a \Borel set $B \subseteq X$ such that:
  \begin{enumerate}
    \item The restriction of $G$ to $B$ is two-regular.
    \item Every connected component of $G$ contains at most one connected
      component of the restriction of $G$ to $B$.
    \item No two points of $B$ of $G$-degree at least three have odd $G$-distance
      from one another.
    \item There is a \Borel perfect matching of $G$ off of $B$.
  \end{enumerate}
  In particular, if there are no injective $G$-rays of $G$-degree two on even indices,
  then the set $B$ is empty, thus $G$ has a \Borel perfect matching.
\end{theorem}

\begin{theoremproof}
  Fix a \Borel set $\calB \subseteq \sets{\N}{X}$ selecting a finite non-empty set of ends
  from every connected component of $G$. By Theorem \ref{Borel:uniformization}, we
  can assume that the set $\calB$ selects exactly $n$ ends from every connected
  component of $G$, for some $n > 0$.
  
  If $n = 1$, then $G$ is generated by the \Borel function $f \from X \to X$ associating to
  each point $x$ its unique $G$-neighbor $y$ for which there is an injective $G$-ray in $\calB$ of
  the form $\sequence{x, y, \ldots}$, in which case Proposition \ref{Borel:function} allows
  us to take $B$ to be the injective part of $f$.
  
  If $n = 2$, then let $A$ denote the set of all points $x$ with two distinct $G$-neighbors
  $y$ for which there are injective $G$-rays in $\calB$ of the form $\sequence{x, y, \ldots}$.
  Proposition \ref{Borel:line} allows us to take $B$ to be the set of $x$ in $A$ for which
  there do not exist $y, z \in A \intersection \connectedcomponent{x}{G}$ of $G$-degree at
  least three having odd $G$-distance from one another.
    
  If $n > 2$, then \cite[Lemma 3.19]{JacksonKechrisLouveau} ensures that the equivalence
  relation generated by $G$ is smooth, in which case Proposition \ref{Borel:smoothmatching}
  allows us to take $B$ to be the empty set.
\end{theoremproof}

\section{Measurable matchings} \label{measurable}

We begin this section with a fact which, despite being quite well known, seems not to have
previously appeared in the form we require.

\begin{proposition}[\Adams, \Jackson-\Kechris-\Louveau] \label{measurable:ends}
  Suppose that $X$ is a \Polish space, $G$ is an acyclic locally countable \Borel
  graph on $X$, and $\mu$ is a \Borel probability measure on $X$ for which the
  equivalence relation generated by $G$ is $\mu$-hyperfinite. Then there are
  $G$-invariant \Borel sets $A, B \subseteq X$ such that:
  \begin{enumerate}
    \item The equivalence relation generated by $G$ is smooth on $A$.
    \item There is a \Borel set selecting a finite non-empty set of ends from every
      connected component of the restriction of $G$ to $B$.
    \item The set $A \union B$ is $\mu$-conull.
  \end{enumerate}
\end{proposition}

\begin{propositionproof}
  This follows from the proof of \cite[Lemma 3.21]{JacksonKechrisLouveau}.
\end{propositionproof}

\begin{remark}
  Although unnecessary for our arguments, it is worth noting that \cite[Lemma 3.19]
  {JacksonKechrisLouveau} allows us to strengthen condition (2) in Proposition
  \ref{measurable:ends} to the existence of a \Borel set selecting one or two ends
  from every connected component of the restriction of $G$ to $B$.
\end{remark}

\begin{remark}
  It is also worth noting that, using a fairly straightforward metamathematical
  argument, Proposition \ref{measurable:ends} can also be established from \cite
  [Lemma 3.21]{JacksonKechrisLouveau} itself, as opposed to its proof. But
  this approach seems rather needlessly roundabout.
\end{remark}

As a corollary, we obtain the following.

\begin{theorem} \label{measurable:Lebesgue}
  Suppose that $X$ is a \Polish space, $G$ is an acyclic locally countable \Borel
  graph on $X$ of degree at least two and with no injective $G$-rays of $G$-degree
  two on even indices, and $\mu$ is a \Borel probability measure on $X$ for which
  the equivalence relation generated by $G$ is $\mu$-hyperfinite. Then there is
  a $\mu$-conull $G$-invariant \Borel set on which $G$ has a \Borel perfect
  matching.
\end{theorem}

\begin{theoremproof}
  Let $A$ and $B$ denote the $G$-invariant \Borel sets whose existence is
  granted by Proposition \ref{measurable:ends}. Proposition \ref{Borel:smoothmatching}
  yields a \Borel perfect matching of $G$ on $A$, and Proposition \ref{Borel:ends} yields
  a \Borel perfect matching of $G$ on $B$, thus there is a \Borel perfect matching of
  $G$ on $A \union B$.
\end{theoremproof}

In the context of \Baire category, we obtain the analogous result for locally finite graphs.

\begin{theorem} \label{measurable:Baire}
  Suppose that $X$ is a \Polish space and $G$ is an acyclic locally finite \Borel
  graph on $X$ of degree at least two and with no injective $G$-rays of $G$-degree
  two on even indices. Then there is a comeager $G$-invariant \Borel set on which
  $G$ has a \Borel perfect matching.
\end{theorem}

\begin{theoremproof}
  Let $\calX$ denote the set of pairs $\pair{S}{T}$ of finite subsets of $X$, where $S
  \subseteq T$ and both are contained in a connected component of $G$. This set
  inherits a standard \Borel structure from $X$. Let $\calG$ denote the graph on $\calX$
  given by
  \begin{equation*}
    \calG = \set{\pair{\pair{S}{T}}{\pair{S'}{T'}} \in \calX \times \calX}[\pair{S}{T} \neq
      \pair{S'}{T'} \mathand T \cap T' \neq \emptyset].
  \end{equation*}
  By \cite[Proposition 3]{ConleyMiller}, there is a \Borel $\N$-coloring $c$ of $\calG$.

  We will now define a decreasing sequence $\sequence{X_s}[s \in \Bairetree]$
  of \Borel subsets of $X$ such that the graph $\restriction{G}{X_s}$ has degree
  at least two and no injective $(\restriction{G}{X_s})$-ray has 
  $(\restriction{G}{X_s})$-degree two on even indices, for all $s \in \Bairetree$.
  We will simultaneously produce an increasing sequence $\sequence{\iota_s}[s
  \in \Bairetree]$ of \Borel matchings of $G$ such that the domain of $\iota_s$ is
  $X \setminus X_s$, for all $s \in \Bairetree$.
  
  Once we have constructed these, for each $p \in \Bairespace$ we will define $X_p
  = \intersection[n \in \N][X_{\restriction{p}{n}}]$ and $\iota_p \from X \setminus
  X_p \to X \setminus X_p$ by $\iota_p(x) = \iota_{\restriction{p}{n}}(x)$, where $n
  \in \N$ is sufficiently large that $x \in X \setminus X_{\restriction{p}{n}}$. As each
  $\iota_p$ is necessarily a \Borel matching of $G$ with domain $X \setminus X_p$,
  it will only remain to show that the details of our construction ensure the existence
  of $p \in \Bairespace$ for which the saturation of $X_p$ with respect to the
  equivalence relation generated by $G$ is meager.

  We begin by setting $X_\emptystring = X$ and $\iota_\emptystring = \emptyfunction$.
  Suppose now that we have already defined $X_s$ and $\iota_s$. Let $\calX_s$ denote
  the set of pairs $\pair{S}{T} \in \calX$ which satisfy the following conditions:
  \begin{enumerate}
    \item The inclusion $\boundary{\restriction{G}{X_s}}{X_s \setminus S} \subseteq T
      \subseteq X_s$ holds.
    \item The graph $\restriction{G}{(X_s \setminus S)}$ has degree at least two.
    \item No injective $(\restriction{G}{(X_s \setminus S)})$-path passing through both a
      point in $\boundary{\restriction{G}{X_s}}{X_s \setminus S}$ and a point in $\boundary
      {\restriction{G}{X_s}}{T}$ has $(\restriction{G}{(X_s \setminus S)})$-degree two on
      even indices.
    \item There is a perfect matching of $\restriction{G}{S}$.
  \end{enumerate}
  We will extend $\iota_s$ by adding perfect matchings of the graphs $\restriction
  {G}{S}$ in a \Borel fashion, using the fact that the sets $T$ provide buffers preventing
  these new matchings from interacting with one another, at least among pairs $\pair{S}
  {T}$ on which our coloring $c$ of $\calG$ is constant.
  
  For each $i \in \N$, define $X_{s \concatenation \sequence{i}} \subseteq
  X_s$ by
  \begin{equation*}
    X_{s \concatenation \sequence{i}} = X_s \setminus \bigcup \set{S}[\exists
      T \ (\pair{S}{T} \in \calX_s \mathand c(S, T) = i)].
  \end{equation*}
  Theorem \ref{Borel:uniformization} ensures that these sets are \Borel.

  \begin{lemma}
    Suppose that $i \in \N$ and $s \in \Bairetree$. Then $\restriction{G}
    {X_{s \concatenation \sequence{i}}}$ has degree at least two.
  \end{lemma}

  \begin{lemmaproof}
    Suppose that $x \in X_{s \concatenation \sequence{i}}$. If $x$ is not in
    $\boundary{\restriction{G}{X_s}}{X_{s \concatenation \sequence{i}}}$,
    then $\degree{\restriction{G}{X_{s \concatenation \sequence{i}}}}{x} = 
    \degree{\restriction{G}{X_s}}{x} \ge 2$. If $x$ is in $\boundary{\restriction{G}
    {X_s}}{X_{s \concatenation \sequence{i}}}$, then there exists $\pair{S}{T} \in
    \calX_s$ such that $c(S, T) = i$ and $x$ is in $\boundary{\restriction{G}{X_s}}
    {X_s \setminus S}$. Condition (1) ensures that $x \in T$, and since $c$ is an
    $\N$-coloring of $\calG$, it follows that $\pair{S}{T}$ is the unique such pair.
    Condition (2) therefore implies that $\degree{\restriction{G}{X_{s \concatenation
    \sequence{i}}}}{x} = \degree{\restriction{G}{(X_s \setminus S)}}{x} \ge 2$.
  \end{lemmaproof}

  \begin{lemma}
    Suppose that $i \in \N$ and $s \in \Bairetree$. Then there is no injective
    $(\restriction{G}{X_{s \concatenation \sequence{i}}})$-ray of $(\restriction{G}
    {X_{s \concatenation \sequence{i}}})$-degree two on even indices.
  \end{lemma}

  \begin{lemmaproof}
    If $\sequence{x_k}[k \in \N]$ is an injective $(\restriction{G}{X_{s \concatenation
    \sequence{i}}})$-ray of $(\restriction{G}{X_{s \concatenation \sequence{i}}})$-degree
    two on even indices, then condition (3) ensures that $x_k \notin \boundary{\restriction
    {G}{X_s}}{X_{s \concatenation \sequence{i}}}$ for all $k \in \N$, thus $\sequence{x_k}
    [k \in \N]$ is a $(\restriction{G}{X_s})$-ray of $(\restriction{G}{X_s})$-degree two on
    even indices, a contradiction.
  \end{lemmaproof}

  Condition (4) ensures that $\iota_s$ extends to a \Borel matching
  $\iota_{s \concatenation \sequence{i}}$ of $G$ with domain $X
  \setminus X_{s \concatenation \sequence{i}}$.

  As noted earlier, it only remains to show that there exists $p \in \Bairespace$
  for which the saturation of $X_p$ with respect to the equivalence relation
  generated by $G$ is meager. We will establish the stronger fact that
  $\forcomeagerlymany p \in \Bairespace \forcomeagerlymany x \in X
  \ \connectedcomponent{x}{G} \intersection X_p = \emptyset$. Note that $\set
  {\pair{p}{x} \in \Bairespace \times X}[x \in X_p]$ is \Borel. By the \Kuratowski-\Ulam
  Theorem (see, for example, \cite[Theorem 8.41]{Kechris}), it is enough to show
  that $\forall x \in X \forcomeagerlymany p \in \Bairespace \ x \notin X_p$. For this,
  it is enough to show that $\forall s \in \Bairetree \forall x \in X_s \exists i \in \N
  \ x \notin X_{s \concatenation \sequence{i}}$.

  Towards this end, suppose that $s \in \Bairetree$ and $x \in X_s$.

  \begin{lemma} \label{lemma:closure}
    There is a finite set $S \subseteq X_s$ such that $x \in S$, $\restriction{G}{(X_s
    \setminus S)}$ has degree at least two, and $\restriction{G}{S}$ has a perfect
    matching.
  \end{lemma}

  \begin{lemmaproof}
    We say that a set $Y \subseteq X$ is \emph{$G$-connected} if the graph
    $\restriction{G}{Y}$ is connected. We will recursively construct increasing 
    sequences $\sequence{\iota_k}[k \in \N]$ of matchings of $G$ and
    $\sequence{S_k}[k \in \N]$ of finite $G$-connected subsets of $X_s$
    containing $x$ such that the domain of $\iota_k$ is $S_k$, for all $k \in \N$.
    We begin by fixing $y \in X_s$ for which $x \mathrel{G} y$, and setting $\iota_0
    = \transposition{x}{y}$ and $S_0 = \set{x, y}$. Given $\iota_k$ and $S_k$,
    observe that for each connected component $C$ of $\restriction{G}{(X_s
    \setminus S_k)}$, there is at most one point $z \in C$ such that $\cardinality
    {C \intersection \verticalsection{G}{z}} = 1$. Let $\iota_{k+1}$ denote the
    minimal extension of $\iota_k$ to an involution which associates every such
    $z$ with the unique element of $C \intersection \verticalsection{G}{z}$, and
    let $S_{k+1}$ denote the domain of $\iota_{k+1}$. This completes the recursive
    construction.

    Set $\iota = \union[k \in \N][\iota_k]$ and $S = \union[k \in \N][S_k]$. Clearly
    $x$ is in $S$, the restriction of $G$ to $X_s \setminus S$ has degree at
    least two, and $\iota$ is a perfect matching of $\restriction{G}{S}$, so it only remains
    to show that $S$ is finite. But if $S$ is infinite, then we can recursively construct an
    injective $(\restriction{G}{X_s})$-ray $\sequence{x_{2k}}[k \in \N]$ with the property that $S
    \intersection \connectedcomponent{x_{2k}}{\restriction{G}{(X_s \setminus S_k)}}$ is
    infinite and $x_{2k + 1}$ is the unique $G$-neighbor of $x_{2k}$ in $X_s \setminus
    S_k$, for all $k \in \N$. But then $\sequence{x_k}[k \in \N]$ has
    $(\restriction{G}{X_s})$-degree two on even indices, a contradiction.
  \end{lemmaproof}

  \begin{lemma} \label{lemma:length}
    There is a finite set $T \subseteq X_s$, with $S \union \boundary{\restriction{G}
    {X_s}}{X_s \setminus S} \subseteq T$, such that no injective $(\restriction{G}{(X_s
    \setminus S)})$-path passing through both a point in $\boundary{\restriction{G}
    {X_s}}{X_s \setminus S}$ and a point in $\boundary{\restriction{G}{X_s}}{T}$ has
    $(\restriction{G}{(X_s \setminus S)})$-degree two on even indices.
  \end{lemma}

  \begin{lemmaproof}
    It is enough to show that for all $z \in X_s \setminus S$, there exists $n \in \N$
    such that there is no injective $(\restriction{G}{(X_s \setminus S)})$-path beginning at $z$,
    having $(\restriction{G}{(X_s \setminus S)})$-degree two on even indices, and along which
    there are $n$ points. Towards this end, observe that if there are arbitrarily long injective
    $(\restriction{G}{(X_s \setminus S)})$-paths of $(\restriction{G}{(X_s \setminus S)})$-degree
    two on even indices beginning at some point $x_0$, then we can recursively
    choose $x_n \notin \set{x_m}[m < n]$ such that there are arbitrarily long injective
    $(\restriction{G}{(X_s \setminus S)})$-paths of $(\restriction{G}{(X_s \setminus S)})$-degree
    two on even indices extending $\sequence{x_k}[k \le n]$, in which case $\sequence{x_k}
    [k \in \N]$ is an injective $(\restriction{G}{(X_s \setminus S)})$-ray of
    $(\restriction{G}{(X_s \setminus S)})$-degree two on even indices, a contradiction.
  \end{lemmaproof}

  As $\pair{S}{T} \in \calX_s$, it follows that $i = c(S, T)$ is as desired.
\end{theoremproof}

We close the paper by noting that the above result fails in the more general locally
countable setting.

\begin{theorem} \label{measurable:example}
  There is a \Polish space $X$ and an acyclic $\aleph_0$-regular \Borel graph $G$ on
  $X$ which does not have a \Borel perfect matching on a comeager \Borel set.
\end{theorem}

\begin{theoremproof}
  We will find \Polish spaces $X$ and $Y$ and a \Borel set $R \subseteq X \times
  Y$, whose horizontal and vertical sections are countably infinite, such that for no
  comeager \Borel set $C \subseteq X$ is there a \Borel injection $\phi \from C \to Y$
  whose graph is contained in $R$. For such $R$, let $G_R$ denote the graph on the
  disjoint union of $X$ and $Y$ in which two points are related if one of them is in $X$,
  one of them is in $Y$, and the corresponding pair is in $R$. As long as we are able
  to simultaneously ensure that $G_R$ is acyclic, it will have the desired properties.
  
  Towards this end, we will recursively define $R_n \subseteq S_n \subseteq
  \Cantorspace[n] \times \Cantorspace[n]$. The sets $R_n$ will provide increasingly
  precise approximations to the set $R$ we seek, whereas the sets $S_n$ will provide
  restrictions on the construction aimed at ruling out the existence of injections whose
  graphs are contained in $R$. We will obtain $R_{n+1}$ and $S_{n+1}$ from the sets
  \begin{equation*}
    R_{n+1}' = \set{\pair{u \concatenation \sequence{i}}{v \concatenation \sequence{i}}}
      [i < 2 \mathand u \mathrel{R_n} v]
  \end{equation*}
  \centerline{and}
  \begin{equation*}
    S_{n+1}' = \set{\pair{u \concatenation \sequence{i}}{v \concatenation \sequence{j}}}
      [i, j < 2 \mathand u \mathrel{R_n} v]
  \end{equation*}
  by either adding a pair to the former or subtracting pairs from the latter, depending on
  whether $n$ is even or odd. Define a function $\projection[0] \from \Cantortree \times
  \Cantortree \to \Cantortree$ by setting $\projection[0](u, v) = u$. In order to ensure that
  the construction can continue, we will proceed in such a fashion that for all $n \in \N$,
  the following conditions hold:
  \begin{enumerate}
    \item $\forall u \in \Cantorspace[n] \exists v \in \Cantorspace[n]
      \ u \mathrel{S_n} v$.
    \item $\forall v \in \Cantorspace[n] \exists u \in \Cantorspace[n]
      \setminus \image{\projection[0]}{R_n} \ u \mathrel{S_n} v$.
  \end{enumerate}
  
  We will describe the exact fashion in which this is accomplished in terms of sequences
  $u_n, v_n \in \Cantorspace[n]$, for $n \in \N$. We can already define the sequences of
  the form $u_{2n}$ and $v_{2n+1}$. In fact, these need only be chosen in such a fashion
  that the corresponding sets $\set{u_{2n}}[n \in \N]$ and $\set{v_{2n+1}}[n \in \N]$
  are \definedterm{dense} in $\Cantortree$, in the sense that
  \begin{equation*}
    \forall t \in \Cantortree \exists m, n \in \N \ t \extendedby u_{2m}, v_{2n+1}.
  \end{equation*}
  The remaining sequences will be chosen during the construction.
  
  We begin by setting $R_0 = \emptysequence$ and $S_0 = \set{\pair{\emptysequence}
  {\emptystring}}$. Suppose now that $n \in \N$ and we have already found $R_{2n}$ and
  $S_{2n}$ satisfying conditions (1) and (2). By the former, there exists $v_{2n} \in
  \Cantorspace[2n]$ such that $u_{2n} \mathrel{S_{2n}} v_{2n}$. It then follows that the sets
  \begin{equation*}
    R_{2n+1} = R_{2n+1}' \union \set{\pair{u_{2n} \concatenation \sequence{0}}{v_{2n}
      \concatenation \sequence{1}}}
  \end{equation*}
  and $S_{2n+1} = S_{2n+1}'$ satisfy conditions (1) and (2) as well. By the latter, there exists
  $u_{2n+1} \in \Cantorspace[2n+1] \setminus \image{\projection[0]}{R_{2n+1}}$ such that
  $u_{2n+1} \mathrel{S_{2n+1}} v_{2n+1}$. It then follows that the sets $R_{2n+2} = R_{2n+2}'$
  and
  \begin{equation*}
    S_{2n+2} = S_{2n+2}' \setminus \set{\pair{u_{2n+1} \concatenation \sequence{0}}{v}}
      [v \neq v_{2n+1} \concatenation \sequence{1}]
  \end{equation*}
  also satisfy conditions (1) and (2). This completes the recursive construction. Define
  $R \subseteq \Cantorspace \times \Cantorspace$ by
  \begin{equation*}
    R = \set{\pair{u \concatenation x}{v \concatenation x}}[n \in \N \mathcomma u \mathrel
      {R_n} v \mathcommaand x \in \Cantorspace].
  \end{equation*}
  
  A simple induction using the definition of $R_n$ reveals that each of the graphs
  $G_{R_n}$ is acyclic, from which it follows that so too is $G_R$.
  
  Another simple induction utilizing the density of $\set{u_{2n}}[n \in \N]$ and $\set{v_{2n+1}}[n
  \in \N]$ along with the definition of $R_n$ and $S_n$ ensures that the sets
  $U_n = \set{x \in \Cantorspace}[\cardinality{\verticalsection{R}{x}} \ge n]$ and $V_n = \set
  {y \in \Cantorspace}[\cardinality{\horizontalsection{R}{y}} \ge n]$ are dense and open.
  Fix homeomorphisms $\phi_n$ of $\Cantorspace$ with the property that $\Ezero = \union[n \in
  \N][\graph{\phi_n}]$, and note that the set $Z = \intersection[m, n \in \N][\preimage{\phi_m}
  {U_n \intersection V_n}]$ is a countable intersection of dense open sets, so the subspace
  topology on $Z$ is \Polish (see, for example, \cite[Theorem 3.11]{Kechris}), and a subset of
  $Z$ is comeager if and only if it is comeager when viewed as a subset of $\Cantorspace$.
  Set $X = Y = Z$, and observe that when viewed as a subset of $X \times Y$, every
  horizontal and vertical section of $R$ is countably infinite.
  
  Suppose, towards a contradiction, that there is a comeager \Borel set $C \subseteq
  \Cantorspace$ for which there is a \Borel injection $\phi \from C \to \Cantorspace$
  whose graph is contained in $R$. As the definitions of $R_n$ and $G$ ensure that
  \begin{equation*} \textstyle
    \graph{\phi} \subseteq \union[n \in \N][{\set{\pair{u_{2n} \concatenation \sequence{0}
      \concatenation x}{v_{2n} \concatenation \sequence{1} \concatenation x}}[x \in \Cantorspace]}],
  \end{equation*}
  there exists $n \in \N$ for which the set of $x \in \Cantorspace$ with the property that
  $\phi(u_{2n} \concatenation \sequence{0} \concatenation x) = v_{2n} \concatenation
  \sequence{1} \concatenation x$ is non-meager. By localization (see, for example,
  \cite[Proposition 8.26]{Kechris}), there exists $r \in \Cantortree$ such that $\phi(u_{2n}
  \concatenation \sequence{0} \concatenation r \concatenation x) = v_{2n} \concatenation
  \sequence{1} \concatenation r \concatenation x$, for comeagerly many $x \in
  \Cantorspace$. Fix $m \ge n$ with $v_{2n} \concatenation \sequence{1}
  \concatenation r \extendedby v_{2m+1}$, noting that $u_{2m+1}$ is incompatible with
  $u_{2n} \concatenation \sequence{0}$ and $\phi$ necessarily sends sequences beginning with
  $u_{2m+1} \concatenation \sequence{0}$ to sequences beginning with $v_{2m+1} \concatenation 
  \sequence{1}$. Again appealing to the above restriction on the graph of $\phi$ imposed
  by the definitions of $R_n$ and $G$, there exists $\ell > m$ for which $u_{2m+1}
  \concatenation \sequence{0} \extendedby u_{2\ell}$ and the set of $x \in \Cantorspace$
  with the property that $\phi(u_{2\ell} \concatenation \sequence{0} \concatenation x) =
  v_{2\ell} \concatenation \sequence{1} \concatenation x$ is non-meager. By one more
  appeal to localization, there exists $s \in \Cantortree$ with the property that $\phi(u_{2\ell}
  \concatenation \sequence{0} \concatenation s \concatenation x) = v_{2\ell} \concatenation 
  \sequence{1} \concatenation s \concatenation x$, for comeagerly many $x \in \Cantorspace$.
  As $u_{2m+1} \concatenation \sequence{0} \extendedby u_{2\ell}$, it follows that $v_{2m+1}
  \concatenation \sequence{1} \extendedby v_{2\ell}$. And since $v_{2n} \concatenation \sequence{1} 
  \concatenation r \extendedby v_{2m + 1}$, it follows that $v_{2n} \concatenation \sequence{1} 
  \concatenation r \extendedby v_{2\ell}$. Fix $t \in \Cantortree$ such that $v_{2\ell} \concatenation
  \sequence{1} \concatenation s = v_{2n} \concatenation \sequence{1} \concatenation r
  \concatenation t$, and observe that
  \begin{align*}
    \phi(u_{2\ell} \concatenation \sequence{0} \concatenation s \concatenation x)
      & = v_{2\ell} \concatenation \sequence{1} \concatenation s \concatenation x \\
      & = v_{2n} \concatenation \sequence{1} \concatenation r \concatenation t
        \concatenation x \\
      & = \phi(u_{2n} \concatenation \sequence{0} \concatenation r \concatenation
        t \concatenation x),
  \end{align*}
  for comeagerly many $x \in \Cantorspace$. As $u_{2m+1}$ is incompatible with
  $u_{2n} \concatenation \sequence{0}$ and extended by $u_{2\ell}$, it follows that
  $u_{2\ell}$ is incompatible with $u_{2n} \concatenation \sequence{0}$, thus
  $u_{2\ell} \concatenation \sequence{0} \concatenation s \concatenation x$
  and $u_{2n} \concatenation \sequence{0} \concatenation r \concatenation t \concatenation
  x$ are distinct sequences with the same image under $\phi$, a contradiction.
\end{theoremproof}

\begin{remark}
  The same idea can, in fact, be used to rule out the existence of a comeager \Borel set
  $C \subseteq \Cantorspace$ for which there is a finite-to-one \Borel function $\phi \from
  C \to Y$ whose graph is contained in $R$.
\end{remark}

\begin{remark}
  It is tempting to try to strengthen the conclusion of Theorem \ref{measurable:example}
  to show that any \Borel matching of $G$ has meager domain, but this is impossible.
  Indeed, any locally countable \Borel graph has a countable \Borel edge coloring \cite
  [Proposition 4.10]{KechrisSoleckiTodorcevic}, and at least one of these colors must
  use a non-meager set of vertices.
\end{remark}

\bibliographystyle{amsalpha}
\bibliography{bibliography}

\providecommand{\bysame}{\leavevmode\hbox to3em{\hrulefill}\thinspace}
\providecommand{\MR}{\relax\ifhmode\unskip\space\fi MR }
\providecommand{\MRhref}[2]{%
  \href{http://www.ams.org/mathscinet-getitem?mr=#1}{#2}
}
\providecommand{\href}[2]{#2}
\begin{thebibliography}{HKL90}

\bibitem[CL]{CsokaLippner}
E.~Cs\'oka and G.~Lippner, \emph{Invariant random matchings in {C}ayley
  graphs}, preprint available at http://arxiv.org/abs/1211.2374.

\bibitem[CM]{ConleyMiller}
C.T. Conley and B.D. Miller, \emph{A bound on measurable chromatic numbers of
  locally finite {B}orel graphs}, to appear in Math. Res. Letters, preprint
  available at http://www.logic.univie.ac.at/benjamin.miller.

\bibitem[DJK94]{DoughertyJacksonKechris}
R.L. Dougherty, S.C. Jackson, and A.S. Kechris, \emph{The structure of
  hyperfinite {B}orel equivalence relations}, Trans. Amer. Math. Soc.
  \textbf{341} (1994), no.~1, 193--225. \MR{1149121 (94c:03066)}

\bibitem[HKL90]{HarringtonKechrisLouveau}
L.A. Harrington, A.S. Kechris, and A.~Louveau, \emph{A {G}limm-{E}ffros
  dichotomy for {B}orel equivalence relations}, J. Amer. Math. Soc. \textbf{3}
  (1990), no.~4, 903--928.

\bibitem[HM09]{HjorthMiller}
G.~Hjorth and B.D. Miller, \emph{Ends of graphed equivalence relations. {II}},
  Israel J. Math. \textbf{169} (2009), 393--415. \MR{2460911 (2009j:03077)}

\bibitem[JKL02]{JacksonKechrisLouveau}
S.~Jackson, A.S. Kechris, and A.~Louveau, \emph{Countable {B}orel equivalence
  relations}, J. Math. Log. \textbf{2} (2002), no.~1, 1--80. \MR{1900547
  (2003f:03066)}

\bibitem[Kec95]{Kechris}
A.S. Kechris, \emph{Classical descriptive set theory}, Graduate Texts in
  Mathematics, vol. 156, Springer-Verlag, New York, 1995. \MR{1321597
  (96e:03057)}

\bibitem[KM]{KechrisMarks}
A.S. Kechris and A.S. Marks, \emph{Descriptive graph combinatorics}, Available
  at http://www.its.caltech.edu/{\twiddle}marks/.

\bibitem[KM04]{KechrisMiller}
A.S. Kechris and B.D. Miller, \emph{Topics in orbit equivalence}, Lecture Notes
  in Mathematics, vol. 1852, Springer-Verlag, Berlin, 2004. \MR{2095154
  (2005f:37010)}

\bibitem[KST99]{KechrisSoleckiTodorcevic}
A.S. Kechris, S.~Solecki, and S.~Todorcevic, \emph{Borel chromatic numbers},
  Adv. Math. \textbf{141} (1999), no.~1, 1--44. \MR{MR1667145 (2000e:03132)}

\bibitem[LN11]{LyonsNazarov}
R.~Lyons and F.~Nazarov, \emph{Perfect matchings as {IID} factors on
  non-amenable groups}, European J. Combin. \textbf{32} (2011), no.~7,
  1115--1125. \MR{2825538 (2012m:05423)}

\bibitem[Mar]{Marks}
A.S. Marks, \emph{A determinacy approach to {B}orel combinatorics}, preprint
  available at http://www.its.caltech.edu/$\sim$marks/.

\bibitem[Mil12]{Miller:Survey}
B.D. Miller, \emph{The graph-theoretic approach to descriptive set theory},
  Bull. Symbolic Logic \textbf{18} (2012), no.~4, 554--575. \MR{3053069}

\bibitem[MU]{MarksUnger}
A.S. Marks and S.~Unger, \emph{Baire measurable paradoxical decompositions via
  matchings}, preprint available at
  http://www.its.caltech.edu/{\twiddle}marks/.

\end{thebibliography}

\end{document}